\documentclass[11pt,leqno]
{amsart}
\usepackage{amsmath,epsfig,graphicx,color}
\usepackage{float}
\definecolor{darkblue}{rgb}{.2, 0.2,.8}
\definecolor{darkgreen}{rgb}{0,0.5,0.3}
\definecolor{darkred}{rgb}{.8, .1,.1}

\newcommand{\clt}{central limit theorem}

\newcommand{\sta}{St\u aric\u a}

\newcommand{\asy}{asymptotic}

\newtheorem{lemma}{Lemma}[section]

\newtheorem{theorem}[lemma]{Theorem}

\newtheorem{proposition}[lemma]{Proposition}
\newtheorem{definition}[lemma]{Definition}
\newtheorem{corollary}[lemma]{Corollary}
\newtheorem{example}[lemma]{Example}
\newtheorem{exercise}[lemma]{Exercise}
\newtheorem{remark}[lemma]{Remark}
\newtheorem{fig}[lemma]{Figure}
\newtheorem{tab}[lemma]{Table}

\newcommand{\bfM}{{\bf M}}

\newcommand{\RV}{{\rm RV}}
\newcommand{\LN}{{\rm LN}}
\newcommand{\WE}{{\rm WE}}
\newcommand{\bth}{\begin{theorem}}
\newcommand{\ethe}{\end{theorem}}

\newcommand{\bre}{\begin{remark}\em }
\newcommand{\ere}{\end{remark}}

\newcommand{\ble}{\begin{lemma}}
\newcommand{\ele}{\end{lemma}}

\newcommand{\pp}{point process}
\newcommand{\bde}{\begin{definition}}
\newcommand{\ede}{\end{definition}}
\newcommand{\bco}{\begin{corollary}}
\newcommand{\eco}{\end{corollary}}

\newcommand{\bpr}{\begin{proposition}}
\newcommand{\epr}{\end{proposition}}

\newcommand{\bexer}{\begin{exercise}}
\newcommand{\eexer}{\end{exercise}}

\newcommand{\bexam}{\begin{example}\rm }
\newcommand{\eexam}{\end{example}}

\newcommand{\bfi}{\begin{fig}}
\newcommand{\efi}{\end{fig}}

\newcommand{\btab}{\begin{tab}}
\newcommand{\etab}{\end{tab}}

\newcommand{\fidi}{finite-dimensional distribution}
\newcommand{\rv}{random variable}

\newcommand{\MDA}{{\rm MDA}}

\newcommand{\PRM}{{\rm PRM}}

\newcommand{\var}{{\rm var}}

\newcommand{\rhs}{right-hand side}

\newcommand{\beao}{\begin{eqnarray*}}
\newcommand{\eeao}{\end{eqnarray*}\noindent}

\newcommand{\beam}{\begin{eqnarray}}
\newcommand{\eeam}{\end{eqnarray}\noindent}

\newcommand{\beqq}{\begin{equation}}
\newcommand{\eeqq}{\end{equation}\noindent}

\newcommand{\bce}{\begin{center}}
\newcommand{\ece}{\end{center}}

\newcommand{\diag}{{\rm diag}}
\newcommand{\barr}{\begin{array}}
\newcommand{\earr}{\end{array}}

\newcommand{\stp}{\stackrel{\P}{\rightarrow}}
\newcommand{\std}{\stackrel{d}{\rightarrow}}

\newcommand{\stv}{\stackrel{v}{\rightarrow}}

\newcommand{\vague}{\stackrel{\lower0.2ex\hbox{$\scriptscriptstyle
                    \it{v} $}}{\rightarrow}}
\newcommand{\weak}{\stackrel{\lower0.2ex\hbox{$\scriptscriptstyle
                    \it{w} $}}{\rightarrow}}
\newcommand{\what}{\stackrel{\lower0.2ex\hbox{$\scriptscriptstyle
                    \it{\hat{w}} $}}{\rightarrow}}

\newcommand{\bdis}{\begin{displaymath}}
\newcommand{\edis}{\end{displaymath}\noindent}

\newcommand{\R}{\mathbb{R}}

\newcommand{\nto}{n\to\infty}

\newcommand{\xto}{x\to\infty}

\newcommand{\ov}{\overline}
\newcommand{\wt}{\widetilde}

\newcommand{\vep}{\varepsilon}

\newcommand{\la}{\lambda}

\newcommand{\regvary}{regularly varying}
\newcommand{\slvary}{slowly varying}
\newcommand{\regvar}{regular variation}

\newcommand{\bbr}{{\mathbb R}}

\newcommand{\con}{convergence}

\newcommand{\evd}{extreme value distribution}

\newcommand{\st}{such that}
\newcommand{\fif}{if and only if}

\newcommand{\fct}{function}

\newcommand{\ds}{distribution}

\newcommand{\cmt}{continuous mapping theorem}
\newcommand{\seq}{sequence}

\newcommand{\ins}{insurance}

\newcommand{\pro}{probabilit}

\newcommand{\ms}{measure}
\newcommand{\mgf}{moment generating function}
\newcommand{\ld}{large deviation}
\newcommand{\bfx}{{\bf x}}
\newcommand{\bfX}{{\bf X}}

\newcommand{\bfy}{{\bf y}}

\newcommand{\bfS}{{\bf S}}

\newcommand{\E }{{\mathbb E}}
\renewcommand{\P }{{\mathbb P}}

\newcommand{\1}{{\mathbf 1}}

\allowdisplaybreaks

\begin{document}
\today
\bibliographystyle{plain}
\title[Gumbel and Fr\'echt convergence of 
the maxima of independent random walks]{ Gumbel and Fr\'echet convergence of 
the maxima of independent random walks}
\thanks{Thomas Mikosch's research is partly support by an Alexander von 
Humboldt Research Award. He takes pleasure in thanking the Faculty of
Mathematics of Ruhruniversit\"at Bochum for hosting him in the period 
December 2018--May 2019.}

\author[T. Mikosch]{Thomas Mikosch}
\address{Department  of Mathematics,
University of Copenhagen,
Universitetsparken 5,
DK-2100 Copenhagen,
Denmark}
\email{mikosch@math.ku.dk}
\author[J. Yslas]{Jorge Yslas}
\address{Department  of Mathematics,
University of Copenhagen,
Universitetsparken 5,
DK-2100 Copenhagen,
Denmark}
\email{jorge.yslas1@gmail.com}

\begin{abstract}
We consider \pp\ \con\ for \seq s of iid random walks. The objective
is to derive \asy\ theory for the largest extremes of these random walks.
We show \con\ of the maximum random walk to the Gumbel or the Fr\'echet
\ds s. The proofs heavily depend on precise \ld\ results for
sums of independent \rv s with a finite \mgf\ or with a subexponential
\ds .
\end{abstract}
\maketitle

\section{Introduction}\setcounter{equation}{0}
Let $(X_i)$ be an iid \seq\ of \rv s with generic element $X$, \ds\ $F$ and right 
tail $\ov F=1-F$.
Define the corresponding partial sum process 
\beao
S_0=0\,,\qquad S_n=X_1+\cdots+X_n\,,\qquad n\ge 1\,.
\eeao
Consider iid copies $(S_{ni})_{i=1,2,\ldots}$ of $S_n$. We also 
introduce an integer \seq\ $(p_n)$ \st\ $p=p_n\to\infty$ as $\nto$.
We are interested in the limiting behavior of the $k$ largest values
among $(S_{ni})_{i=1,\ldots,p}$, in particular in the possible limit laws of the 
maximum $\max_{i=1,\ldots,p} S_{ni}$.
More generally, writing $\vep_x$ for Dirac measure at $x$, we are interested
in the limiting behavior of the \pp es 
\beam\label{eq:conv}
N_p=\sum_{i=1}^p \vep_{c_n^{-1} (S_{ni}-d_n)} \std N\,,\qquad \nto\,,
\eeam
for suitable constants $c_n>0$ and $d_n\in\bbr$ toward a 
Poisson random \ms\ $N$ with Radon mean \ms\ $\mu$ (we write $\PRM(\mu)$).
\par
 Our main motivation for this work comes 
from random matrix theory, in particular when dealing with sample covariance
matrices. Their entries are dependent random walks. However, 
in various situations the theory can be modified in such a way
that it suffices to study independent random walks. We refer to 
Section~\ref{subsec:related} for a discussion.
\par
Relation \eqref{eq:conv} is equivalent to the following limit relations
for the tails
\beao
p_n\,\P\big(c_n^{-1}(S_n-d_n)\in (a,b]\big) \to \mu(a,b]\,,
\eeao
for any $a<b$ provided that $\mu(a,b]<\infty$; see Resnick \cite{resnick:2007}, Theorem 5.3. These conditions 
involve precise \ld\ \pro ies for the random walk $(S_n)$; 
in Section~\ref{sec:ld} we provide some results which are relevant 
in this context. 
\par
We distinguish between two types of precise \ld\ results:
\begin{itemize}
\item 
normal approximation
\item
subexponential approximation
\end{itemize}
The normal approximation can be understood as extension of
the \clt\ for $(S_n/\sqrt{n})$ toward increasing intervals.
This approximation causes the maxima of $(S_{ni}/\sqrt{n})$ to behave like the maxima of an iid normal \seq , i.e., these maxima converge in \ds\ to the Gumbel \ds . 
This is in contrast to the subexponential approximation which requires
that $F$ is a so-called {\em subexponential} \ds ; 
see Section~\ref{subsec:subexp}. In particular, 
$\ov F$ is heavy-tailed in the sense that the
\mgf\ does not exist. This fact implies that $\P(S_n>x_n)\sim n\,\ov F(x_n)$
for sufficiently fast increasing \seq s $x_n\to\infty$.
Hence $n\,\ov F(x_n)$ dominates
$\P(S_n>x_n)$ at sufficiently high levels $x_n$ and,  as in
limit theory for the maxima of an iid \seq ,  $\ov F$ 
determines the type of the limit \ds\ of the maxima of $(S_{ni})$ as well
as the normalizing and centering constants. In this case we also
assume that $F$ belongs to the maximum domain of attraction (MDA) of the 
Gumbel or Fr\'echet \ds s, and we borrow the known normalizing and centering
constants from these MDA. Thus, in the case of the MDA of the Gumbel \ds\ 
the maxima of $(S_{ni})$ may converge to the Gumbel \ds\ due to two distinct
mechanisms: the normal approximation at medium-high thresholds or the 
subexponential approximation at high-level thresholds. In the case of the 
MDA of the Fr\'echet \ds\ two distinct approximations are possible:
Gumbel approximation at medium-high thresholds and Fr\'echet approximation at 
high-level thresholds provided the \ds\ has finite second moment. If this condition is not satisfied only the Fr\'echet approximation is possible.
\par
The paper is organized as follows. In Section~\ref{sec:prelim}
we introduce the necessary notions for this paper: subexponential and 
\regvary\ \ds s (Section~\ref{subsec:subexp}), maximum domain of attraction
and relevant \ds s in it (Section~\ref{subsec:mda}), point process \con\
of triangular arrays toward Poisson random \ms s (Section~\ref{subsec:pp}),
precise \ld s (Section~\ref{subsec:ld}). Due to the importance of 
the latter topic we devote Section~\ref{sec:ld} to it and collect
some of the known precise \ld\ results in the case when 
the \mgf\ is finite in some neighborhood of the origin and for
subexponential \ds s. The main results of this paper
are formulated in Section~\ref{sec:main}. Based on the \ld\ results of
Section~\ref{sec:ld} we give sufficient conditions for the \pp\ 
\con\ relation \eqref{eq:conv} to hold and we clarify which rates of
growth are possible for $p_n\to\infty$. In particular, we consider the 
case when $p_n$ in \eqref{eq:conv} is replaced by $k_n=[n/r_n]$ for some integer \seq\
$r_n\to\infty$ and $n$ is replaced by $r_n$. This means that 
we are interested in \eqref{eq:conv} when $S_{ni}=S_{r_ni}-S_{r_n(i-1)}$,
$i=1,\ldots,k_n$, are iid block sums. We also discuss extensions 
of these results to stationary \regvary\ \seq s 
(Section~\ref{subsec:stat}) 
and iid multivariate \regvary\ \seq s (Section~\ref{eq:subsecmultrv}).

\section{Preliminaries I}\label{sec:prelim}\setcounter{equation}{0}
\subsection{Subexponential and \regvary\ \ds s}\label{subsec:subexp}
We are interested in the 
class $\mathcal S$ of {\em subexponential \ds s} $F$, i.e.,  
it is a \ds\ supported on $[0,\infty)$
\st\ for any $n\ge 2$,
\beao
\P(S_n>x)\sim n\,\ov F(x)\,,\qquad \xto\,.
\eeao 
For an encyclopedic treatment of subexponential \ds s, see Foss et al.
\cite{foss:korshunov:zachary:2013}.
In \ins\ mathematics,  ${\mathcal S}$ is considered
a natural class of heavy-tailed \ds s. In particular, $F$ does not
have a finite \mgf ; see Embrechts et al. \cite{embrechts:kluppelberg:mikosch:1997}, Lemma~1.3.5.
\par
The \regvary\ \ds s are another class of  heavy-tailed \ds s
supported on $\bbr$. We say that $X$ and its \ds\ $F$ are {\em \regvary\
with index} $\alpha>0$ if there are a \slvary\ \fct\ $L$ and constants
$p_\pm$ \st\ $p_++p_-=1$ and 
\beam\label{eq:regvar}
F(-x)\sim p_-\,x^{-\alpha}\,L(x)\qquad\mbox{and}\qquad \ov F(x)\sim 
  p_+\,x^{-\alpha}\,L(x)\,,\qquad \xto\,.
\eeam
A non-negative \regvary\ $X$ is subexponential; see \cite{embrechts:kluppelberg:mikosch:1997}, Corollary~1.3.2.
\subsection{Maximum domains of attraction}\label{subsec:mda}
We call a non-degenerate \ds\ $H$ an {\em \evd} if there exist 
constants $c_n>0$ and  $d_n\in\bbr$, $n\ge 1$, \st\ 
the maxima $M_n=\max(X_1,\ldots,X_n)$
satisfy the limit relation
\beam\label{eq:mda}
c_n^{-1} (M_n-d_n) \std Y\sim H\,,\qquad\nto\,. 
\eeam
In the context of this paper we will deal
with two standard \evd s: the Fr\'echet \ds\ $\Phi_\alpha(x)=\exp(-x^{-\alpha})$,
$x>0$, and the Gumbel \ds\ $\Lambda(x)= \exp(-\exp(-x))$, $x\in\bbr$.
As a matter of fact, the third type of \evd\ -- the Weibull \ds\ -- cannot 
appear since \eqref{eq:mda} is only possible for $X$ with finite right endpoint
but a random walk is not bounded from above by a constant.  
We say that the \ds\ $F$ of $X$ is in the {\em maximum domain of attraction} 
of the \evd\ $H$ $(F\in\MDA(H))$.
\bexam A \ds\ $F\in\MDA(\Phi_\alpha)$ for some $\alpha>0$ \fif\ 
\beao
\ov F(x)= \dfrac{L(x)}{x^\alpha}\,,\qquad x>0\,;
\eeao
see \cite{embrechts:kluppelberg:mikosch:1997}, Section~3.3.1.
Then
\beao
c_n^{-1} M_n\std Y\sim \Phi_\alpha\,,\qquad \nto\,,
\eeao
where $(c_n)$ can be chosen \st\ $n\,\P(X>c_n)\to 1$.
\eexam
\bexam\label{exam:mdagumbel} A \ds\ $F$ with infinite right endpoint obeys $F\in\MDA(\Lambda)$  \fif\ there exists a positive \fct\ $a(x)$ with derivative $a'(x)\to 0$ as $\xto$
\st
\beao
\lim_{u\to\infty}\dfrac{\ov F(u+ a(u)\, x)}{\ov F(u)}= {\rm e}^{-x}\,,\qquad x\in\bbr\,;
\eeao
see \cite{embrechts:kluppelberg:mikosch:1997}, Section~3.3.3.
Then
\beao
c_n^{-1} (M_n-d_n)\std Y\sim \Lambda\,,\qquad \nto\,,
\eeao
where $(d_n)$ can be chosen \st\ $n\,\P(X>d_n)\to 1$ and $c_n=a(d_n)$.
\par
The standard normal \ds\ $\Phi\in \MDA(\Lambda)$ and satisfies
\beam\label{eq:23}
c_n^{-1} (M_n-d_n)\std Y\sim \Lambda\,,\qquad\nto\,,
\eeam
where $c_n=1/d_n$ and 
\beam\label{eq:constn}
d_n=\sqrt{2\log n}-\frac{\log\log n+\log 4\pi}{2(2\log n)^{1/2}}\,.
\eeam 
Since $d_n\sim \sqrt{2\log n}$ we can replace $c_n$ in \eqref{eq:23} 
by $1/\sqrt{2\log n}$ while $d_n$ cannot be replaced by $\sqrt{2\log n}$.
\par
The standard lognormal \ds\ (i.e., $X=\exp(Y)$ for a standard normal \rv\ $Y$)
is also in $\MDA(\Lambda)$. In particular,
one can choose 
\beam\label{eq:lognormal}
c_n=d_n/\sqrt{2\log n}\quad\mbox{and}\quad
d_n=\exp\Big(\sqrt{2\log n} - \frac{\log\log n+\log 4\pi}{2(2\log n)^{1/2}}\Big)\,;
\eeam
see \cite{embrechts:kluppelberg:mikosch:1997}, p. 156.
\par
The standard Weibull \ds\ has tail $\ov F(x)= \exp(-x^{-\tau})$, $x>0$, $\tau>0$.
We consider a \ds\ $F$ on $(0,\infty)$ with a Weibull-type tail $\ov F(x)\sim 
c\,x^\beta \exp(-\la x^\tau)$ for constants $c,\beta,\la,\tau>0$.
Then $F\in \MDA(\Lambda)$ and one can choose
\beam\label{eq:weibullcndn}
c_n=(\la \tau)^{-1} s_n^{1/\tau-1}\quad\mbox{and}\quad
d_n=s_n^{1/\tau}+\dfrac 1\tau s_n^{1/\tau-1} \Big(\dfrac{\beta}{\la \tau} \log s_n + \dfrac{\log c}{\la}\Big)\,,
\eeam 
where $s_n=\la^{-1}\log n$; see \cite{embrechts:kluppelberg:mikosch:1997},  p.~155.
\eexam
\subsection{Point process \con\ of independent triangular arrays}\label{subsec:pp}
For further use we will need the following \pp\ limit result
(Resnick \cite{resnick:2007}, Theorem~5.3).
\bpr\label{prop:conv}
Let $(X_{ni})_{n=1,2,\ldots;i=1,2,\ldots}$ be a triangular array of row-wise iid random 
elements on some state space $E\subset\bbr^d $ equipped with the Borel $\sigma$-field $\mathcal E$.
Let $\mu$ be a Radon \ms\ on $\mathcal E $.
Then 
\beao
\wt N_p= \sum_{i=1}^p\vep_{X_{ni}}\std N\,,\qquad \nto\,,
\eeao
holds
for some $\PRM(\mu)$ $N$ \fif\ 
\beam\label{eq:vague}
p\,\P(X_{n1}\in \cdot)\stv \mu(\cdot)\,,\qquad \nto\,,
\eeam
where $\stv$ denotes vague \con\ on $E$.
\epr
\subsection{Large deviations}\label{subsec:ld}
Our main goal is to prove the \pp\ \con\ \eqref{eq:conv}
for iid \seq s $(S_{ni})$ of partial sum processes ($\bbr$- or $\bbr^d$-valued), properly normalized and 
centered. It follows from Proposition~\ref{prop:conv}
that this means to prove relations of the type
\beao
p\,\P\big(c_n^{-1}(S_{n}-d_n) \in (a,b])\to \mu(a,b]\quad\mbox{or}\quad
p\,\P\big(c_n^{-1}(S_{n}-d_n) >a)\to \mu(a,\infty)\,,
\eeao
provided $\mu(a,b]+\mu(a,\infty)<\infty$. Since $p=p_n\to\infty$ this means
that $\P\big(c_n^{-1}(S_{n}-d_n) >a)\to 0$ as $\nto$. We will refer 
to these vanishing \pro ies as {\em \ld\ \pro ies}. 
In Section~\ref{sec:ld}
we consider some of the well-known precise \ld\ results 
in heavy- and light-tail situations.

\section{Preliminaries II: precise \ld s}\label{sec:ld}\setcounter{equation}{0}
In this section we collect some precise \ld\ results in the
light- and heavy-tailed cases.
\subsection{Large deviations with normal approximation}
 We assume $\E[X]=0$, $\var(X)=1$ and write
$\Phi$ for the standard normal \ds .
We start with a classical result
when $X$ has finite exponential moments. 
\bth[Petrov's theorem \cite{petrov:1972}, Theorem 1 in Chapter VIII] \label{thm:petrov} 
Assume that the \mgf\ $\E[\exp(h\,X)]$ is finite
in some neighborhood of the origin. Then the following 
tail bound holds for $0\le x=o(\sqrt{n})$:
\beao
\dfrac{\P(S_n/\sqrt{n}>x)}{\ov \Phi(x)}= \exp\Big(\dfrac{x^3}{\sqrt{n}} 
\la\big(\dfrac{x}{\sqrt{n}}\big)\Big) \Big[1+ O\big(\dfrac{x+1}{\sqrt{n}}\big)\Big]\,,\quad \nto\,.
\eeao
where  $\la(t)= \sum_{k=0}^\infty a_kt^k$ is the Cram\'er series whose coefficients $a_k$
depend  on the cumulants of $X$, and $\la(t)$ converges for sufficiently small
values $|t|$.
\ethe
Under the conditions of Theorem~\ref{thm:petrov},  uniformly  for $x=o(n^{1/6})$,
\beam\label{eq:ss}
\dfrac{\P(S_n/\sqrt{n}>x)}{\ov \Phi(x)}\to 1\,,\qquad \nto\,.
\eeam
Theorem 7 in Chapter VIII of Petrov \cite{petrov:1972} considers the 
situation of Theorem~\ref{thm:petrov} under the additional assumption that
the cumulants of order $k=3,\ldots,r+2$ of $X$ vanish for some positive integer $r$. Then the coefficients
$a_0,\ldots,a_{r-1}$ in the series $\la(t)$ vanish, and it is not difficult
to see that \eqref{eq:ss} holds uniformly for $0\le x=o\big(n^{(r+1)/(2(r+3))}\big)$. 
\par
In \cite{petrov:1972}, Section VIII.3, one also finds necessary and sufficient 
conditions for \eqref{eq:ss} to hold in certain intervals.
The following result was proved by S.V. Nagaev \cite{nagaev:1965}
for $x\in (0,\sqrt{(s/2-1)\log n})$ and improved by R. Michel \cite{michel:1974}
for $x\in (0,\sqrt{(s-2)\log n})$.
The statement of the proposition is sharp under the given moment condition; see Theorem~\ref{thm:nagaev} below.
\bpr\label{prop:sd}
Assume that $\E[|X|^s]<\infty$
for some $s>2$. Then \eqref{eq:ss} holds uniformly for $0\le x\le 
\sqrt{(s-2)\log n}$.
\epr
\subsection{Large deviations with normal/subexponential approximations}\label{subsec:ld2}
Cline and Hsing \cite{cline:hsing:1998} (in an unpublished article)
discovered that the subexponential class ${\mathcal S}$ of
\ds s exhibits  a completely different kind of \ld\ behavior: 

\bpr[Cline and Hsing \cite{cline:hsing:1998}]\label{prop:clinehsing}
 We consider a \ds\ $F$ on $(0,\infty)$ with infinite right endpoint.
Then the following statements hold.
\begin{enumerate}
\item[\rm 1.]
$F\in {\mathcal S}$ \fif\ 
\beam\label{eq:longtailed}
\lim_{\xto}\dfrac{\ov F(x+y)}{\ov F(x)}=1\,,\quad\mbox{for any real $y$,}
\eeam
and there exists a \seq\ $\gamma_n\to\infty$ \st
\beam\label{eq:febr18a}
\lim_{\nto} \sup_{x>\gamma_n}\dfrac{\P(S_n>x)}{n\,\ov F(x)}\le 1\,.
\eeam
\item[\rm 2.]
If $F\in{\mathcal S}$ then
there exists a \seq\ $\gamma_n\to\infty$ \st\
\beam\label{eq:febr18b}
\lim_{\nto} \sup_{x>\gamma_n}\Big|\dfrac{\P(S_n>x)}{n\,\ov F(x)}-1\Big|=0\,.
\eeam
\end{enumerate}
\epr
\bre\label{rem:longtailed} If $F$ satisfies \eqref{eq:longtailed} we say that $F$ is 
{\em long-tailed}, we write $F\in {\mathcal L}$.
It is well known that $F\in {\mathcal S}$ implies $F\in \mathcal{L}$; see
Embrechts et al. \cite{embrechts:kluppelberg:mikosch:1997}, 
Lemma 1.3.5 on p.~41. The converse is not true.
\ere
\par
Proposition~\ref{prop:clinehsing} shows that the subexponential class
is the one for which heavy-tail \ld s are reasonable to study. 
Given that we know that $F$ is long-tailed, $F$ is subexponential \fif\ 
a uniform \ld\ relation of the type \eqref{eq:febr18b} holds.

Subexponential and normal approximations to large deviation 
\pro ies were studied in detail in various papers. Among them,
\ld s for iid \regvary\ \rv s are perhaps studied best.
S.V. Nagaev \cite{nagaev:1979} formulated a seminal result
about the \ld s of a random walk $(S_n)$ in the case of
\regvary\ $X$ with finite variance. He dedicated this theorem to his
brother A.V. Nagaev who had started this line of research
in the 1960s; see for example \cite{nagaev:1969a,nagaev:1969}.

\bth[Nagaev's theorem \cite{nagaev:1969a,nagaev:1979}]\label{thm:nagaev}
Consider an iid \seq\ $(X_i)$ of \rv s with $\E[X]=0$, $\var(X)=1$
and $\E[|X|^{2+\delta}]<\infty$ for some $\delta>0$. 
Assume that $\ov F(x)=  x^{-\alpha}\,L(x)$, $x>0$, for some $\alpha>2$ 
and a \slvary\ \fct\ $L$. Then for $x\ge \sqrt{n}$ as $\nto$,
\beao
\P(S_n>x)=\ov \Phi(x/\sqrt{n})\,(1+o(1))+ n\,\ov F(x)\,(1+o(1))\,.
\eeao  
In particular, if $X$ satisfies \eqref{eq:regvar} with
constants $p_\pm$, 
then for any positive constant $c_1< \alpha-2$
\beam\label{eq:apprnorma}
\sup_{1<x/\sqrt{n}< \sqrt{c_1\,\log n}}\Big|
\dfrac{\P(\pm S_n>x)}{\ov \Phi\big(x/\sqrt{n}\big)}-1\Big|\to 0\,,\qquad \nto\,,
\eeam  
and for any constant $c_2>\alpha-2$,
\beao
\sup_{x/\sqrt{n}>\sqrt{c_2\,\log n}}\Big|
\dfrac{\P(\pm S_n>x)}{n\,\P(|X|>x)}-p_{\pm}\Big|\to 0\,,\qquad\nto\,.
\eeao
\ethe
\bre If $X$ is \regvary\ with index $\alpha$, $\E[|X|^{s}]$ is finite (infinite) for   $s<\alpha$ $(s>\alpha)$. Therefore the normal approximation
\eqref{eq:apprnorma} is in agreement with Proposition~\ref{prop:sd}.
\ere

In the infinite variance \regvary\ 
case this result is complemented by an analogous statement. It can be found 
in Cline and Hsing \cite{cline:hsing:1998}, Denisov et al.
\cite{denisov:dieker:shneer:2008}. 
 \bth\label{thm:heyde}
Consider an iid \seq\ $(X_i)$ of \regvary\ \rv s with 
index $\alpha\in (0,2]$ satisfying \eqref{eq:regvar}. Assume $\E[X]=0$ if
this expectation is finite.
Choose $(a_n)$ \st\
\beao
n\,\P(|X|>a_n) + \dfrac{n}{a_n^2}\E[X^2\,\1(|X|\le a_n)]=1\,,\qquad n=1,2,\ldots,
\eeao
and $(\gamma_n)$ \st\ $\gamma_n/a_n\to\infty$ as $\nto$.
For $\alpha=2$, also assume for sufficiently small $\delta>0$,
\beao
\lim_{\nto}\sup_{x>\gamma_n} \dfrac{n}{x^2}\,\dfrac{\E[X^2\,\1(|X|\le x)]}
{[n\,\P(|X|>x)]^\delta}=0\,.
\eeao
Choose $(d_n)$ \st
\beam\label{eq:dn}
d_n= \left\{\barr{ll}0\,,&\quad \alpha\in (0,1)\cup(1,2]\,,\\[2mm]
n\;\E [X\,\1(|X|\le a_n)]\,,&\quad\alpha=1\,.
\earr\right.
\eeam 
Then the following \ld\ result holds:
\beao
\sup_{x>\gamma_n}\Big|\dfrac{\P\big(\pm(S_n-d_n)>x\big)}{n\,\P(|X|>x)}-p_\pm \Big|\to 0\,,\qquad\nto\,.
\eeao  
\ethe
\bre The normalization $(a_n)$ is chosen \st\ $a_n^{-1}(S_n-d_n)\std Y_\alpha$
for an $\alpha$-stable \rv\ $Y_\alpha$, $\alpha\in (0,2]$.
Therefore $\gamma_n^{-1}(S_n-d_n)\stp 0 $.
In the case $\alpha<2$, in view of Karamata's theorem (see Bingham et al.
\cite{bingham:goldie:teugels:1987}), it is possible to  choose $(a_n)$
according as $n\,\P(|X|>a_n)\to 1$. The case $\alpha=2$ is delicate:
in this case $\var(X)$ can be finite or infinite. In the former case,
$(a_n)$ is proportional to $\sqrt{n}$, in the latter case 
$(a_n/\sqrt{n})$ is a \slvary\ \seq ; see Feller \cite{feller:1971} or
Ibragimov and Linnik \cite{ibragimov:linnik:1971}, Section II.6. 
\ere 
Normal and subexponential approximations to \ld\ \pro ies also
exist for subexponential \ds s that have all moments finite. 
Early on, this was observed by A.V. Nagaev \cite{nagaev:1969a,nagaev:1969,nagaev:1977}.
Rozovskii \cite{rozovski:1993} did not use the name of subexponential
\ds , but the conditions on the tails of the \ds s he introduced are
``close'' to subexponentiality; he also allowed for \ds s $F$ supported 
on the whole real line. In particular, A.V. Nagaev and Rozovskii discovered 
that, in general, the $x$-regions where the normal and subexponential
approximations hold are separated from each other. To make this precise,
we call two \seq s $(\xi_n)$ and $(\psi_n)$ {\em separating \seq s for
the normal and subexponential approximations to \ld\ \pro ies} if
for an iid \seq\ $(X_i)$ with variance 1,
\beao
&&\sup_{x< \xi_n}\Big|\dfrac{\P(S_n-\E[S_n]>x)}{\ov \Phi(x/\sqrt{n})}-1\Big|\to 0\,,\\
&&\sup_{x> \psi_n}\Big|\dfrac{\P(S_n-\E[S_n]>x)}{n\,\P(X>x)}-1\Big|\to 0\,,\qquad \nto\,.
\eeao
A.V. Nagaev and  Rozovskii gave conditions under 
which $(\psi_n)$ and $(\xi_n)$ cannot have the same \asy\ order; i.e.,
one necessarily has $\psi_n/\xi_n\to \infty$.
In particular, in the $x$-region $(\xi_n,\psi_n)$ neither the normal
nor the subexponential approximation holds; Rozovskii \cite{rozovski:1993}
also provided \ld\ approximations for $\P(S_n>x)$ for these regions involving 
$\ov \Phi(x/\sqrt{n})$ and a truncated Cram\'er series. Explicit expressions
for $(\psi_n)$ and $(\xi_n)$ are in general hard to get. We focus on two
classes of subexponential \ds s where the separating \seq s are known.
\begin{itemize}
\item
{\em Lognormal-type tails}, we write  $F\in \LN(\gamma)$:
for some constants $\beta,\xi\in \bbr$, $\gamma>1$ and  $\lambda,c>0$,
\beao
\ov F(x)\sim c\,x^{\beta}\,(\log x)^\xi\,\exp\big(-\lambda\,(\log x)^\gamma\big)\,,\qquad x\to\infty\,.
\eeao
In the notation  $\LN(\gamma)$ we suppress the dependence on 
$\beta,\xi,\lambda,c$.
\item
{\em Weibull-type tails}, we write $F\in \WE(\tau)$: for some 
$\beta \in \bbr$, $\tau\in (0,1)$,  $\lambda,c>0$.
\beao
\ov F(x)\sim c\,x^{\beta}\,\exp\big(-\lambda\, x^{\tau}\big)\;,\qquad x\to\infty\,.
\eeao
In the notation $\WE(\tau)$ we suppress the 
dependence on $\beta,\lambda,c$.
\end{itemize}
\par
The name ``Weibull-type tail'' is motivated by the fact
that the Weibull \ds\ $F$ with shape parameter $\tau\in (0,1)$ belongs to 
$\WE(\tau)$. Indeed, in this case $\ov F(x)=\exp(-\lambda x^\tau)$, $x>0$,
for positive parameters $\la$. Similarly, the lognormal \ds\ $F$ belongs to
 $\LN(2)$. This is easily seen by an application of Mill's ratio: for a standard normal \rv\ $Y$,
\beao
\ov F(x)= \P(Y>\log x)\sim \dfrac{\exp\big(- (\log x)^2/2\big)}{\sqrt{2\pi}\,\log x}\,,\qquad \xto\,.
\eeao

\begin{table}[H]
\begin{center}
\begin{tabular}{|l|l|l|}\hline
\;$F\in$\;&\; $\xi_n$\; & \;$\psi_n$\;\\ \hline
\;$\RV(\alpha)$\,,\;$\alpha>2$\; &\;$\big((\alpha-2)n\,\log n\big)^{1/2}$\;&
\;$\big((\alpha-2)n\,\log n\big)^{1/2}$\\[2mm]
\;$\LN(\gamma)$\,,\;$1<\gamma<2$\; &\;$(n\,(\log n)^{\gamma})^{1/2}$\;
&\;$(n\,(\log n) ^{\gamma})^{1/2}$\\[2mm]
\;$\LN(\gamma)$\,,\;$\gamma \ge 2$\; &\;$(n\,(\log n)^{\gamma})^{1/2} /\wt h_n $\;
&\;$n^{1/2} (\log n)^{\gamma-1}\,h_n$\\[2mm]
\;$\WE(\tau)$\,,\;$0<\tau\le 0.5$\; &\;$n^{1/(2-\tau)}/ \wt h_n$\;
&\;$n^{1/(2-2\tau)}\, h_n$\\[2mm]
\;$\WE(\tau)$\,,\;$0.5<\tau<1$\; &\;$n^{2/3}/ \wt h_n$\;
&\;$n^{1/(2-2\tau)}\,h_n$\\
\hline
\end{tabular}
\caption{Separating \seq s $(\xi_n)$ and $(\psi_n)$ for the normal
and subexponential approximations of $\P(S_n-\mathbb{E}[S_n]>x)$. We also assume 
$\var(X)=1$. Here $(h_n), (\wt h_n)$  are any \seq s  converging to infinity.
For completeness, we also include the \regvary\ class $\RV(\alpha)$.
The table is taken from Mikosch and Nagaev~\cite{mikosch:nagaev:1998}.}
\label{tab:1}
\end{center}
\end{table}

These classes of \ds s have rather distinct tail behavior. 
It follows from the theory in  
Embrechts et al. \cite{embrechts:kluppelberg:mikosch:1997}, 
Sections 1.3 and 1.4, that membership of $F$ in $\RV(\alpha)$, $\LN(\gamma)$ or $\WE(\tau)$ implies $F\in\mathcal S$.
The case $\WE(\tau)$, $0<\tau<1$, was already considered by A.V. Nagaev \cite{nagaev:1969,nagaev:1977}.
\par
For the heaviest tails when $F\in \LN(\gamma)$, $1<\gamma< 2$ 
one can still choose $\xi_n=\psi_n$.  This means that one threshold
\seq\ separates the normal and subexponential approximations to the
right tail $\P(S_n-\E[S_n]>x)$.
Rozovskii \cite{rozovski:1993} discovered that 
the classes $\LN(\gamma)$, $\gamma \ge 2$ ,
and $\LN(\gamma)$, $1< \gamma< 2$ have rather distinct \ld\ properties.
In the case $\gamma\ge 2$ one cannot choose 
$(\xi_n)$ and $(\psi_n)$ the same. 
The class \LN($\gamma$) with $1<\gamma < 2$ satisfies the conditions of Theorem 3b in Rozovskii \cite{rozovski:1993} which implies that
\beao
	\P(S_n-\E[S_n]>x ) = \big[ \overline{\Phi}(x/\sqrt{n} )\1 \big(  x<\gamma_n \big)  + n \overline{F}(x) \1 \big( x>\gamma_n \big) \big](1+o(1))
\eeao
uniformly for $x$, where $\gamma_{n} = \big(\lambda2^{-\gamma+1} \big)^{1/2} n^{1/2}(\log n)^{\gamma/2}$.
For $\gamma=2$ the conditions of Theorem 3a in \cite{rozovski:1993} 
are satisfied: with $g(x) =  \lambda (\log x)^{2} - ({\beta}+2) \log x -{\xi}  \log \left( \log x \right) - \log c $ and as $\nto$,
\beao
	 \P(S_n-\E[S_n]>x) = \big[ \overline{\Phi}(x/\sqrt{n} )\,\1(x<\gamma_n)+
	n \overline{F}(x) {\rm e}^{\frac{n(g'(x))^{2}}{2}}\1(x>\gamma_n)\big](1+o(1)) \,.
\eeao
Direct calculation shows that $\P(S_n-\E[S_n]>\gamma_n)\sim \exp(\la)\,n\,\ov F(\gamma_n)$ while, uniformly for $x > \gamma_n h_{n} $, $h_{n} \to \infty $, we have that 
$\P(S_n-\E[S_n]>x ) \sim n \overline{F}(x)$.
\par
It is interesting to observe that all but one class 
of subexponential \ds s considered
in Table~\ref{tab:1} have the property that $c\,n \in (\psi_n,\infty)$ for any $c>0$. The exception is $\WE(\tau)$ for $\tau\in (0.5,1)$. This fact
turns the investigation of the tail \pro ies $\P(S_n-\E[S_n]>c\,n)$ into a 
complicated technical problem.
The exponential ($\WE(1)$) and superexponential ($\WE(\tau)$), $\tau> 1$,
classes do not contain subexponential \ds s. 
The corresponding partial sums exhibit the light-tailed
\ld\ behavior of Petrov's Theorem~\ref{thm:petrov}.
As a historical remark, Linnik \cite{linnik:1961a}
and S.V. Nagaev \cite{nagaev:1965} determined lower separating 
\seq s $(\xi_n)$ for the normal approximation to the tails $\P(S_n-\E[S_n]>x)$
under the assumption that $\ov F$ is dominated by the tail of a  
regular subexponential \ds\ from the table. 
\par
Denisov et al. 
\cite{denisov:dieker:shneer:2008} and Cline and Hsing \cite{cline:hsing:1998}
considered a unified approach to subexponential \ld\ approximations
for general subexponential and related \ds s. In particular, they identified separating \seq s $(\psi_n)$
for the subexponential approximation of the tails $\P(S_n-\E[S_n]>x)$
for general subexponential \ds s.  Denisov et al. 
\cite{denisov:dieker:shneer:2008} also considered local versions, i.e.,
approximations to the tails $\P(S_n\in [x,x+T])$ for $T>0$ as $\xto$.

\section{Main results}\label{sec:main}\setcounter{equation}{0}
\subsection{Gumbel convergence via
normal approximations to \ld\ \pro ies for small $x$.}
We assume that  $\E[X]=0$ and $\var(X)=1$ and the \ld\ approximation to the standard normal \ds\ $\Phi$ holds: for some 
$\gamma_n\to\infty$,
\beam\label{eq:equiv1}
\sup_{0\le x< \gamma_n}\big|\dfrac{\P(S_n/\sqrt{n}>x)}{\ov \Phi(x)}-1\Big|\to 0\,,\qquad \nto\,.
\eeam
\par
We recall that  $\Phi\in \MDA(\Lambda)$ and \eqref{eq:23} holds.
An analogous relation holds for the maxima of iid random walks 
$S_{n1}/\sqrt{n},\ldots, S_{np}/\sqrt{n}$ as follows from the next
result.
\bth\label{thm:1} Assume that \eqref{eq:equiv1} is satisfied for some $\gamma_n\to\infty$. Then 
\beam\label{eq:2u}
p\,\P\Big( \dfrac{S_n}{\sqrt{n}} >d_p+x/d_p\Big) 
\to {\rm e}^{-x}\,,\qquad \nto\,,\qquad x\in\bbr\,, 
\eeam
holds for any integer \seq\ $p_n\to\infty$
\st\ $p_n< \exp(\gamma_n^2/2)$ and  $(d_p)$ is defined in \eqref{eq:constn}. 
Moreover, for the considered $(p_n)$, \eqref{eq:2u} is equivalent to either of the following 
limit relations:
\begin{enumerate}
\item
For $\Gamma_i=E_1+\cdots +E_i$ and an iid standard exponential \seq\ $(E_i)$
the following \pp\ \con\ holds on the state space $\bbr$
\beam\label{eq:22}
N_p=\sum_{i=1}^p \vep_{d_p\,\big(\frac{S_{ni}}{\sqrt{n}} -d_p\big)}\std 
N=\sum_{i=1}^\infty \vep_{-\log \Gamma_i}\,,
\eeam
where $N$ is $\PRM(-\log \Lambda)$ on $\bbr$.
\item Gumbel \con\ of the maximum random walk
\beao
d_p\max_{i=1,\ldots,p} \big(S_{ni}/\sqrt{n}-d_p\big)\std Y\sim
\Lambda\,,\qquad \nto\,.
\eeao
\end{enumerate}
\ethe
\begin{proof} 
In view of Proposition~\ref{prop:conv} it suffices for $N_p\std N$ 
to show that
\beao
p\,\P\Big( d_{p} \big(\dfrac{S_n}{\sqrt{n}}-d_p\big) >x\Big) =
p\,\P\Big( \dfrac{S_n}{\sqrt{n}} >d_p+x/d_p\Big) \sim p\,\ov \Phi(d_p+x/d_p) 
\to {\rm e}^{-x}\,,\qquad x\in\bbr\,.
\eeao
But this follows from \eqref{eq:equiv1} and the definition of $(d_p)$ 
if we assume that
$d_p+x/d_p<\gamma_n$, i.e., $p_n< \exp(\gamma_n^2/2)$
\st\ $p_n\to \infty$. 
\par
If $N_p\std N$ a continuous mapping argument implies that 
\beao
\P\big(N_p(x,\infty)=0\big)&=& 
\P\Big(\max_{i=1,\ldots,p} d_p\big(S_{ni}/\sqrt{n}-d_p\big)\le x\Big)\\
&\to&\P\big(N(x,\infty)=0\big)= \Lambda(x)\,,\qquad x\in\bbr\,,\qquad\nto\,.
\eeao
On the other hand, for $x\in\bbr$ as $\nto$,
\beao
\P\Big(\max_{i=1,\ldots,p} d_p\big(S_{ni}/\sqrt{n}-d_p\big)\le x\Big)
&=& \Big(1- \dfrac {p\,\P\big(d_p (S_{n1}/\sqrt{n}-d_p\big)>x\big)}{p}\Big)^p \\
&\to& \exp\big(-{\rm e}^{-x}\big)\,,
\eeao
\fif\ \eqref{eq:2u} holds.
\end{proof}
\bre
If one replaces in  \eqref{eq:22} the quantities 
$(S_{ni}/\sqrt{n})_{i=1,\ldots,p}$ by iid standard normal
\rv s  then this limit relation remains valid.
This means that, under \eqref{eq:equiv1}, e.g. under 
the assumption of a finite  \mgf\ in some neighborhood
of the origin (see Section~\ref{sec:ld}),  
the \clt\ makes the tails of $(S_{ni}/\sqrt{n})_{i=1,\ldots,p}$
almost indistinguishable from those of the standard normal \ds .
This is in stark contrast to subexponential \ds s where the 
characteristics of $\ov F(x)$ show up in the tail $\P(S_{ni}/\sqrt{n}>x)$
for large values of $x$.
\ere
\subsubsection{The extreme values of iid random walks.}\label{rem:1a}
Write 
\beao
S_{n,(p)}\le \cdots \le S_{n,(1)}
\eeao
for the
ordered values of $S_{n1},\ldots,S_{np}$
The following result is immediate from Theorem~\ref{thm:1}.


\bco\label{cor:analog}
Assume that the conditions of Theorem~\ref{thm:1} hold.
Then
\beam\label{eq:ko}
\sqrt{2\log p}\,\Big(\dfrac{S_{n,(1)}}{\sqrt{n}} -d_p,\ldots ,\dfrac{S_{n,(k)}}{\sqrt{n}} -d_p \Big)\std \big(- \log \Gamma_1,\ldots,-\log \Gamma_k\big) \,,\qquad
\nto\,.
\eeam
 Moreover, if \eqref{eq:equiv1} also holds for the sequence $(-X_{i})$, then we have 
 
 Moreover, if there is $\gamma_{n} \to \infty $  such that   $\sup_{0\le x< \gamma_n}\big|{\P( \pm S_n/\sqrt{n}>x)}/{\ov \Phi(x)}-1\Big|\to 0$ as $\nto$,  then we have
\beam\label{eq:k1}
\lefteqn{
\P\Big(\max_{i=1,\ldots,p} d_p\big(S_{ni}/\sqrt{n}-d_p\big)\le x\,,\min_{i=1,\ldots,p} d_p\big(S_{ni}/\sqrt{n}+d_p\big)\le y\Big)
}\nonumber\\&\to &\Lambda(x) \big(1-\Lambda(-y)\big)\,,\qquad x,y\in\bbr\,,\qquad \nto\,.
\eeam
\eco
\begin{proof}
We observe that $d_p/\sqrt{2\log p}\to 1$. Then  \eqref{eq:22}
and the \cmt\ imply that \eqref{eq:ko} holds
for any fixed $k\ge 1$.
\par
We observe that
\beao\lefteqn{
\P\Big(\max_{i=1,\ldots,p} d_p\big(S_{ni}/\sqrt{n}-d_p\big)\le x\,,\min_{i=1,\ldots,p} d_p\big(S_{ni}/\sqrt{n}+d_p\big)\le y\Big)
}\\&=&\P\Big(\max_{i=1,\ldots,p} d_p\big(S_{ni}/\sqrt{n}-d_p\big)\le x\Big)\\
&&-\P\Big(\max_{i=1,\ldots,p} d_p\big(S_{ni}/\sqrt{n}-d_p\big)\le x\,,\min_{i=1,\ldots,p} d_p\big(S_{ni}/\sqrt{n}+d_p\big)>y\Big)\\
&=&P_1(x,y)-P_2(x,y)\,.
\eeao
Of course, $P_1(x,y)\to \Lambda(x)$. 
 On the other hand, 
\beao
P_2(x,y)
&=&\P\Big(\bigcap_{i=1}^p \big\{ S_{ni}/\sqrt{n}\le d_p+x/d_p\,,
S_{ni}/\sqrt{n}>-d_p+y/d_p
\big\}\Big)\\
&=&\Big( \P\Big( -d_p+y/d_p< S_{n1}/\sqrt{n} \le 
d_p+x/d_p
\Big) \Big)^{p}\\
&=& \exp \Big( p \, \log \Big( 1- \P\Big(  S_{n1}/\sqrt{n} > d_p+x/d_p \Big) -\P\Big( S_{n1}/\sqrt{n} \le -d_p+y/d_p \Big) \Big) \Big)\\
&\to & \exp\big(-({\rm e}^{-x}+{\rm e}^{y})\big)=\Lambda(x)\Lambda(-y)\,.
\eeao
The last step follows from a Taylor expansion of the logarithm and Theorem~\ref{thm:1}. This proves \eqref{eq:k1}.
\end{proof}

\subsubsection{Examples.}\label{subsubsec:examples}
In this section we verify the assumptions of Theorem~\ref{thm:1} for various classes of \ds s $F$. We always assume $\E[X]=0$ and $\var(X)=1$.
\bexam\label{exam:osi}
Assume the existence of the \mgf\ of $X$ in some neighborhood of the 
origin. Petrov's Theorem~\ref{thm:petrov} ensures 
\eqref{eq:22} for $p\le \exp(o(n^{1/3}))$.
\eexam
\bexam\label{exam:moment}
Assume $\E[|X|^s]<\infty$ for some $s>2$. Proposition~\ref{prop:sd} ensures that
\eqref{eq:22} for
  $p\le n^{(s-2)/2}$. 
\eexam
\bexam\label{exam:regvar}
Assume that $X$ is \regvary\ with index $\alpha>2$.
Then we can apply Nagaev's Theorem~\ref{thm:nagaev} with 
$\gamma_n= \sqrt{c\,\log n}$ for any $c<\alpha-2$ and \eqref{eq:22} holds for  $p\le  n^{c/2}$. This is in agreement 
with Example~\ref{exam:moment}.
\eexam
\bexam\label{exam:LN}
Assume that $X$ has a \ds\ in $\LN(\gamma)$ for some $\gamma>1$. From Table~\ref{tab:1}, $\gamma_n=o((\log n )^{\gamma/2}$,  and \eqref{eq:22} holds  for
$p\le \exp( o((\log n )^{\gamma}))$ 
\eexam
\bexam\label{exam:WE}
Assume that $F\in \WE(\tau)$, $0<\tau<1$.
Table~\ref{tab:1} yields $\gamma_n=o(n^{\tau/(2(2-\tau))})$  for $\tau\le 0.5$, 
hence $p\le \exp(o(n^{\tau/(2-\tau)}))$, and  for $\tau\in (0.5,1)$,
$\gamma_n=o(n^{1/6})$ and $p\le \exp(o(n^{1/3}))$.
\eexam
We summarize these examples in Table~\ref{tab:2}.

\begin{table}[H]
\begin{center}
\begin{tabular}{|l|l|}\hline
Example No & Upper bound for $p$\\
\hline
\ref{exam:osi} Petrov case & $\exp(o(n^{1/3}))$\\[2mm]
\ref{exam:moment} $\mathbb{E}[|X|^s]<\infty$, $s>2$& $n^{(s-2)/2}$\\[2mm]
\ref{exam:regvar} ${\rm RV}(\alpha)$, $\alpha>2$, $c<\alpha-2$& $ n^{c/2}$\\[2mm]
\ref{exam:LN} $\LN(\gamma)$, $\gamma>1$ &$\exp( o((\log n )^{\gamma}))$ \\[2mm]
\ref{exam:WE}  $\WE(\tau)$, $\tau\le 0.5$&$\exp(o(n^{\tau/(2-\tau)}))$ \\[2mm]
\ref{exam:WE} $\WE(\tau)$, $\tau\in (0.5,1)$&   $\exp(o(n^{1/3}))$  \\
\hline 
\end{tabular}
\caption{Upper bounds for $p$.}
\label{tab:2}
\end{center}
\end{table}

\subsubsection{The extremes of the blocks of a random walk.}\label{exam:sumsn}\rm We consider a random walk $S_n$ with iid step sizes
$X_i$ with $\E[X]=0$ and $\var(X)=1$, and with \ds\ $F$, and any integer \seq\ $r_n\to \infty$ \st\
$k_n=[n/r_n]\to \infty$ as $\nto$. Set $S_{ni}= S_{r_n i} -S_{r_n(i-1)}$, i.e.,
this is the sum of the $i$th block $X_{r_n(i-1)+1},\ldots,X_{r_n i}$.
Then we are in the setting of Theorem~\ref{thm:1} if we replace 
$p_n$ by $k_n$ and $n$ by $r_n$. We are interested in the following 
result for the \pp\ of the block sums of $S_n$ with length $r_n$ (see \eqref{eq:22}) 
\beam\label{eq:NP2}
N_{k_n}=\sum_{i=1}^{k_n} \vep_{d_{k_n}\,\big(\frac{ S_{r_n i} -S_{r_n(i-1)}}{\sqrt{r_n}} -d_{k_n}\big)}\std 
N=\sum_{i=1}^\infty \vep_{-\log \Gamma_i}\,.
\eeam
This means we are looking for $(r_n)$
\st\ $n/r_n<\exp(\gamma_{r_n}^2/2)$. This amounts to the following conditions
on $(r_n)$ in Table~\ref{tab:3}:

\begin{table}[H]
\begin{center}
\begin{tabular}{|l|l|}\hline 
Example No & Lower bounds for $(r_n)$\\
\hline 
\ref{exam:osi} Petrov case & $r_n/(\log n)^3\to \infty$\\[2mm]
\ref{exam:regvar} $\RV(\alpha)$, $\alpha>2$& $r_n> n^{2/(\alpha-\vep)}$ any $\vep\in (0,\alpha-2)$ \\[2mm]
\ref{exam:moment} $\mathbb{E}[|X|^s]<\infty$, $s>2$& $r_n>n^{2/s}$\\[2mm]
\ref{exam:LN} $\LN(\gamma)$, $\gamma>1$ & $r_n/ \exp((2 \log n)^{1/\gamma}) \to \infty $\\[2mm]
\ref{exam:WE} $\WE(\tau)$,  $\tau\le 0.5$& $r_n /(\log n)^{(2-\tau)/\tau}\to\infty$\\[2mm]
\ref{exam:WE} $\WE(\tau)$,  $\tau\in (0.5,1)$&$r_n/(\log n)^{3}\to\infty$\\ 
\hline 
\end {tabular} 
\caption{Lower bounds on the block size $r_n$}
\label{tab:3}
\ece
\end{table}

This table shows convincingly that, the heavier the tails, the larger we have to choose the block length $r_n$. Otherwise, the normal approximation does
not function sufficiently well simultaneously 
for the block sums $S_{r_n i} -S_{r_n(i-1)}$, $i=1,\ldots,k_n$. In particular, 
in the \regvary\ case we always need that $r_n$ grows polynomially.
\par
Notice that we have from \eqref{eq:NP2} in particular
\beao
\dfrac{d_{k_n}}{\sqrt{r_n}} \,\max_{i=1,\ldots,k_n}  \Big(S_{r_n i} -S_{r_n(i-1)}-\sqrt{r_n}\,d_{k_n}\Big) \std -\log \Gamma_1\sim \Lambda\,,\qquad\nto\,.
\eeao
The normalization $d_{k_n}/\sqrt{r_n}$
is \asy\ to $\sqrt{(2 \log k_n)/r_n}$.

\subsection{Gumbel convergence via the subexponential approximation to \ld\ \pro ies for
very large $x$}
In this section we will exploit the subexponential approximation
to \ld\ \pro ies for subexponential \ds s $F$, i.e.,
\beam\label{eq:subexpld}
\sup_{x>\gamma_n} \Big| \dfrac{\P(S_n-\E[S_n]>x)}{n\,\P(X>x)}- 1\Big|\to 0\,,
\eeam
and we will also assume that $F\in\MDA(\Lambda)$; see Example~\ref{exam:mdagumbel} for the corresponding MDA conditions and the definition of the
centering constants $(d_n)$ and the normalizing constants $(c_n)$.
Then, in particular,
$X$ has all moments finite. In this case, the Gumbel approximation of
the \pp\ of the $(S_{ni})$ is also possible.
\bth\label{thm:x} Assume that $F\in\MDA(\Lambda)\cap \mathcal S$, 
the subexponential approximation \eqref{eq:subexpld} holds
and for sufficiently large $n$ and an integer \seq\ $p_n\to\infty$,
\beam\label{eq:larger}
d_{np}+x\,c_{np}>\gamma_n \,, \qquad  \mbox{for any $x<0$\,,}
\eeam
 where $(d_{np})$ and $(c_{np})$ are the subsequences of $(d_{n})$ and $(c_{n})$, respectively, evaluated at $n p$. Then
\beam\label{eq:analog2}
p\,\P\big(S_n-\E[S_n]>d_{np}+ x\, c_{np}\big) \to {\rm e}^{-x}\,,\qquad x\in\bbr\,,\qquad \nto\,,
\eeam
holds. Moreover, \eqref{eq:analog2} is equivalent to either of the 
following  limit relations:
\begin{enumerate}
\item
Point process \con\ to a Poisson process on the state space $\bbr$
\beam\label{eq:gumbel2}
N_p= \sum_{i=1}^p\vep_{c_{np}^{-1} (S_{ni}-\E[S_n]-d_{np})}\std N\,,\qquad \nto\,,
\eeam
where $N\sim\PRM(-\log \Lambda)$; see Theorem~\ref{thm:1}.
\item
Gumbel \con\ of the maximum random walk
\beam\label{eq:analog3}
\max_{i=1,\ldots,p} c_{np}^{-1}\big((S_{ni}-\E[S_n])-d_{np}\big)\std Y\sim \Lambda\,,\qquad \nto\,.
\eeam
\end{enumerate}
\ethe 
\begin{proof} If $d_{np}+x\,c_{np}>\gamma_n$ for every $x<0$ it holds for 
$x\in\bbr$. Therefore \eqref{eq:subexpld} applies.
Since  $F\in\MDA(\Lambda)\cap {\mathcal S}$ and by definition 
of $(c_n)$ and $(d_n)$ we have
\beao
p\,\P\big(S_n-\E[S_n]>d_{np}+   x\, c_{np}\big) \sim 
p\,n\,\P(X>d_{np} +x\,c_{np})\to {\rm e}^{-x}\,,\quad x\in\bbr\,,\qquad \nto\,,
\eeao
proving \eqref{eq:analog2}.
Proposition~\ref{prop:conv} yields the equivalence of 
\eqref{eq:gumbel2} and  \eqref{eq:analog2}.
The equivalence of \eqref{eq:gumbel2} and \eqref{eq:analog3}
follows from a standard argument.
\end{proof}
\bre\label{rem:delta}
Since $a(x)$ defined in Example~\ref{exam:mdagumbel} has density $a'(x)\to 0$ as $\xto$ we have $a(x)/x\to 0$.
On the other hand, $c_n=a(d_n)$ and $d_n\to\infty$ since $F\in\mathcal S$.
Therefore for any $x>0$, 
\beao
d_{np}+x\,c_{np}= d_{np}\Big(1+ x\,\dfrac{a(d_{np})}{d_{np}}\Big)
\sim d_{np}\,.
\eeao 
Hence \eqref{eq:larger} holds if $d_{np}\ge (1+\delta)\gamma_n$ for 
any small $\delta>0$ and large $n$.
\ere

\subsubsection{The extreme values of iid random walks.}
Relation \eqref{eq:gumbel2} and a continuous mapping argument
imply the following analog of Corollary~\ref{cor:analog}.
We use the same notation as in  Section~\ref{rem:1a}. One can follow
the lines of  the proof of Corollary~\ref{cor:analog}.

\bco Assume the conditions of Theorem~\ref{thm:x}.
Then the following relation holds for $k\ge 1$,
\beao
c_{np}^{-1}\Big(S_{n,(1)}-\E[S_n] -d_{np},\ldots ,
S_{n,(k)}-\E[S_n] -d_{np} \Big)\std \big(- \log \Gamma_1,\ldots,-\log \Gamma_k\big) 
\eeao
as $\nto$.
\eco
\subsubsection{Examples.}
Theorem~\ref{thm:x}  applies to $F\in\LN(\gamma)$, $\gamma>1$, and 
$F\in \WE(\tau)$, $0<\tau<1$; see the discussion in Section~\ref{subsec:ld2}.
However, the calculation of the constants $(c_n)$ and $(d_n)$ is rather
complicated for these classes of subexponential \ds s. For illustration of the 
theory we restrict 
ourselves to two parametric classes of \ds s where these constants are known.
\bexam
We assume that $X$ has a standard lognormal \ds . From \eqref{eq:lognormal},
Table~\ref{tab:1} and Remark~\ref{rem:delta} we conclude that we need to verify the condition $
\exp\big(\sqrt{2\log (np)}\big)\ge h_n\sqrt{n} \log n$
for a \seq\ $(h_n)$ increasing to infinity arbitrarily slowly.
Calculation shows that it suffices to 
choose $p_n\to\infty$ \st\ 
$p> \exp\big((\log n)^2\big)$.
\eexam

\bexam We assume that $X$ has a Weibull \ds\ with tail 
$\ov F(x)=\exp({-x^\tau})$ for some $\tau\in (0,1)$. From \eqref{eq:weibullcndn} we conclude that $d_{np}\sim (\log np)^{1/\tau}$. In view of
Remark~\ref{rem:delta} and Table~\ref{tab:1} it suffices to verify that
$
 (\log np)^{1/\tau}\ge h_n \, n^{1/(2-2\tau)}
$ 
for a \seq\ $h_n\to\infty$ arbitrarily slowly. It holds if
$
p> n^{-1}\,\exp\big(\big(h_n n^{1/(2-2\tau)}\big)^\tau\big)
$.
\eexam

\subsubsection{The extremes of the blocks of a random walk.}
We appeal to the notation in Section~\ref{exam:sumsn}.
We are in the setting of Theorem~\ref{thm:x} if we replace 
$p_n$ by $k_n$ and $n$ by $r_n$. We are interested in the following 
result for the \pp\ of the block sums of $S_n$ with length $r_n$ (see \eqref{eq:gumbel2}) 
\beam\label{eq:NP2a}
N_{k_n}=\sum_{i=1}^{k_n} \vep_{c_{n}^{-1}\,\big(S_{r_n i} -S_{r_n(i-1)} -\E[S_{r_n}]-d_{n}\big)}\std 
N=\sum_{i=1}^\infty \vep_{-\log \Gamma_i}\,.
\eeam
We need to verify condition \eqref{eq:larger} which turns into
$d_n+c_n\,x>\gamma_{r_n}$. In view of Remark~\ref{rem:delta} it suffices to
prove that $d_n>h_n\gamma_{r_n}$ for a \seq\ $h_n\to\infty$ arbitrarily slowly; see Table~\ref{tab:1} for some $\gamma_{n}$-values. 
\par
We start with a standard lognormal \ds ; see \eqref{eq:lognormal} for the
corresponding $(c_n)$ and $(d_n)$. In particular, we need to verify
\beao
d_n=\exp\Big(\sqrt{2\log n} - \frac{\log\log n+\log 4\pi}{2(2\log n)^{1/2}}\Big)
\ge h_n \sqrt{r_n} \log r_n\,.
\eeao
A sufficient condition is $\exp(2\sqrt{2\log n})>\wt h_n r_n$
for a \seq\ $\wt h_n\to\infty$ arbitrarily slowly. We observe that the left-hand expression is a \slvary\ \fct .
\par
Next we consider a standard Weibull \ds\ for $\tau\in (0,1)$. The constants
$(c_n)$ and $(d_n)$ are given in \eqref{eq:weibullcndn}. In particular,
we need to verify 
\beao
d_n\sim ( \log n)^{1/\tau}> h_n r_n^{1/(2-2\tau)}\,.
\eeao 
This holds if $(\log n)^{2(1-\tau)/\tau} h_n ^{-2(1-\tau)}>r_n$. Again,
this is a strong restriction on the growth of $(r_n)$ and is in contrast
to the \regvary\ case where polynomial growth of $(r_n)$ is possible; 
see Section~\ref{subsubsec:blcokfrechet}.

\subsection{Fr\'echet \con\ via the subexponential approximations to
\ld\ \pro ies for large $x$}
In this section we assume that $X$ is \regvary\ with index $\alpha>0$
in the sense of \eqref{eq:regvar}.
Throughout we choose a  normalizing \seq\ $(a_n)$ \st\
$n\,\P(|X|>a_n)\to 1$ as $\nto$.
The following result is an analog of Theorems~\ref{thm:1} and \ref{thm:x}.

\bth\label{thm:2}
Assume that $X$ is \regvary\ with index $\alpha>0$ and $\E[X]=0$ if the expectation is finite. Choose a \seq\ $(d_n)$ \st\ 
\beao
d_n= \left\{\barr{ll}0\,,&\quad \alpha\in (0,1)\cup(1,\infty)\,,\\[2mm]
n\;\E [X\,\1(|X|\le a_n)]\,,&\quad\alpha=1\,,
\earr\right.
\eeao 
We assume that $p_n\to\infty$ is an integer \seq\ which
satisfies the additional conditions
\beam\label{eq:nag}\left\{\barr{ll}
a_{np}\ge \sqrt{(\alpha-2+\delta) n \log n}
\mbox{ for some small $\delta>0$}& \mbox{if $\alpha>2$\,,}\\
\lim_{\nto}\sup_{x>a_{np}} p^\delta \dfrac{n}{x^2}\,\E[X^2\,\1(|X|\le x)]=0
\mbox{ for some small $\delta>0$}&\mbox{if $\alpha=2$\,.}
\earr\right.
\eeam 
Then the following limit relation
\beam\label{eq:analog4}
p\, \P(\pm a_{np}^{-1} (S_n-d_n)>x)\to p_\pm x^{-\alpha}\,,\qquad x>0\,,\qquad\nto\,,
\eeam
holds. Moreover, \eqref{eq:analog4} is equivalent to 
\beam\label{eq:ppcon}
N_p=\sum_{i=1}^p \vep_{a_{np}^{-1}(S_{ni}-d_n)}\std N=\sum_{i=1}^\infty\vep_{q_i\,\Gamma_i^{-1/\alpha}}\,,  
\eeam
where $(\Gamma_i)$ is defined in Theorem~\ref{thm:1} and $(q_i)$ is an iid
\seq\ of Bernoulli variables with \ds\ $\P(q_i=\pm 1)=p_\pm$ independent of 
$(\Gamma_i)$.
\ethe
\begin{proof} We start by verifying \eqref{eq:analog4}.
Assume $\alpha<2$. Then for 
any \seq\ $p_n\to\infty$, $a_{np}/a_n\to \infty$. Therefore  Theorem~\ref{thm:heyde} and the definition of $(a_{np})$ yield
\beao
p\, \P(\pm a_{np}^{-1} (S_n-d_n)>x)\sim p\,n\,\P(\pm X>a_{np}\,x)\sim p_\pm x^{-\alpha}\,,\qquad \nto\,.
\eeao
If $\alpha>2$ the same result holds in view of Theorem~\ref{thm:nagaev} 
since we assume condition \eqref{eq:nag}. If $\alpha=2$ we can again
apply Theorem~\ref{thm:heyde} with $\gamma_n=a_{np}$ and use \eqref{eq:nag}.
\par
We notice that the limit \pp\ $N$ is $\PRM(\mu_\alpha)$
with intensity
\beam\label{eq:mualpha}
\mu_\alpha(dx)= |x|^{-\alpha-1}\big(p_+ \1(x>0)+p_- \1(x<0)\big)\,dx\,.
\eeam
An appeal to Proposition~\ref{prop:conv} shows that \eqref{eq:analog4}
and \eqref{eq:ppcon} are equivalent.

\end{proof}
\bre Assume $\alpha>2$.
Since $a_{np}=(np)^{1/\alpha}\ell(np)$ for a 
\slvary\ \fct\ $\ell$ and $\ell(x)\ge x^{-\gamma/\alpha}$ for any small
$\gamma>0$ and sufficiently large $x$, \eqref{eq:nag} holds
if $p\ge n^{(\alpha/2)-1+\gamma'}$ for any choice of $\gamma'>0$.
Assume $\alpha=2$ and $\var(X)<\infty$. Then $a_{np}\sim c\,\sqrt{np}$ 
and \eqref{eq:nag} is satisfied for any \seq\ $p_n\to\infty$ and $\delta<1$.
If $\var(X)=\infty$, $a_{np}= (np)^{1/2}\ell(np)$ for a \slvary\ \fct\ $\ell$
and $\E[X^2(|X|\le x)]$ is an increasing \slvary\ \fct . Using
Karamata bounds for \slvary\ \fct s, we conclude that \eqref{eq:nag}
holds if $p/n^{\gamma}\to \infty$ for any small $\gamma>0$.
\ere
\subsubsection{The extreme values of iid random walks.}
For simplicity, we assume $d_n=0$.
Write $N_p^+$ for the restriction of $N_p$ to the state space $(0,\infty)$
and $S^+_{n,(1)}$ for the maximum of 
$(S_{n1})_+,\ldots,(S_{np})_+$. We also write $\xi=\min\{i\ge 1: q_i=1\}$
and assume that $\xi$ is independent of $(\Gamma_i)$.
Then  \eqref{eq:ppcon}
and the \cmt\ imply that 
\beam
\P\big(N_p^+(x,\infty)=0\big)&=& 
\P\big(a_{np}^{-1}S^+_{n,(1)}\le x \big)\nonumber\\
&\std & \P\big(\Gamma_\xi^{-1/\alpha}\le x\big)=\Phi_\alpha^{p_+}(x)\,.\label{eq:feb3}
\eeam
Moreover, we have joint \con\ of minima and maxima.
\bco
Assume the conditions of Theorem~\ref{thm:2} and $d_n=0$.
Then
\beao
\lim_{\nto}\P\Big(0<a_{np}^{-1}\max_{i=1,\ldots,p} S_{ni}\le x\,,
-y<a_{np}^{-1}\min_{i=1,\ldots,p} S_{ni}  \Big)= \Phi_\alpha^{p_+}(x)\Phi_\alpha^{p_-}(y)\,,\quad x,y>0\,.
\eeao
\eco
\begin{proof}
We have 
\beao
\P\Big(a_{np}^{-1}\max_{i=1,\ldots,p} S_{ni}\le x\,, 
-y<a_{np}^{-1}\min_{i=1,\ldots,p} S_{ni}\Big)&=&
\P\big( N_p\big( (x,\infty)\cup (-\infty,-y]\big)=0\big)\\
&\to &\P\big( N\big( (x,\infty)\cup (-\infty,-y]\big)=0\big)\\
&=& \exp\big(-( p_+x^{-\alpha}+ p_- y^{-\alpha})\big)\\
&=&\Phi_\alpha^{p_+}(x)\Phi_\alpha^{p_-}(y)
\,,\qquad\nto\,. 
\eeao
\end{proof}
\subsubsection{The extremes of the blocks of a random walk.}\label{subsubsec:blcokfrechet}
We appeal to the notation of Section~\ref{exam:sumsn} and apply 
Theorem~\ref{thm:2} in the case when $n$ is replaced by some integer-sequence
$r_n\to\infty$ \st\ $k_n=[n/r_n]\to\infty$ and $p_n$ is replaced by $k_n$.
We also assume for simplicity that $d_n=0$.
Observing that $a_{np}$ turns into $a_{r_n k_n}\sim a_n$, 
\eqref{eq:ppcon} turns into
\beao
N_{k_n}=\sum_{i=1}^{k_n} \vep_{a_{n}^{-1}(S_{r_n\,i}-S_{r_n(i-1)})}\std N=\sum_{i=1}^\infty\vep_{q_i\,\Gamma_i^{-1/\alpha}}\,,  \qquad \nto\,.
\eeao
\par
For simplicity, we assume $\alpha\ne 2$. If $\alpha<2$ no further restrictions on $(r_n)$ are required. If $\alpha>2$ we have the additional 
growth condition $a_{n}>\sqrt{(\alpha-2+\delta)r_n\,\log r_n}$ for sufficiently large $n$. Since $a_n=n^{1/\alpha}\ell(n)$ for some \slvary\ \fct\ $\ell$ this
amounts to showing that $n^{2/\alpha}\ell^2(n)/ (\alpha-2+\delta)>r_n\log r_n$.
Since any \slvary\ \fct\ satisfies $\ell(n)\ge n^{-\vep}$ for any $\vep>0$ and $n\ge n_0(\vep)$ we get the following sufficient condition on the growth of 
$(r_n)$: for any sufficiently small $\vep>0$, $n^{2/\alpha-\vep}>r_n$. 
This condition ensures that $(r_n)$ is significantly smaller than $n$, and the larger $\alpha$ the more stringent this condition becomes.
\par
An appeal to  \eqref{eq:feb3} yields in particular
\beao
\P\Big(a_{n}^{-1}\max_{i=1,\ldots,k_n}  (S_{r_n\,i}- S_{r_n\,(i-1)})_+\le x \Big)
&\std & \P\big(\Gamma_\xi^{-1/\alpha}\le x\big)=\Phi_\alpha^{p_+}(x)\,,\\
\P\Big(a_{n}^{-1}\max_{i=1,\ldots,k_n}  |S_{r_n\,i}- S_{r_n\,(i-1)}|\le x \Big)
&\std & \P\big(\Gamma_1^{-1/\alpha}\le x\big)=\Phi_\alpha(x)\,,\qquad \nto\,.
\eeao

\subsubsection{Extension to a stationary \regvary\ sequence.}\label{subsec:stat}
In view of classical theory (e.g. Feller~\cite{feller:1971})
 $X$ is \regvary\ with index $\alpha \in (0,2)$ \fif\ 
$a_n^{-1}(S_n-d_n)\std \xi_\alpha$ for an $\alpha$-stable \rv\ $\xi_\alpha$
where one can choose $(a_n)$ \st\ $n\,\P(|X|>a_n)\to 1$ and 
$(d_n)$ as in \eqref{eq:dn}. For the sake of argument we also assume 
$d_n=0$; this is a restriction only in the case $\alpha=1$.
\par
If $(r_n)$ is any integer \seq\ \st\ $r_n\to\infty$
and $k_n=[n/r_n]\to 0$ then 
\beam\label{eq:sum}
a_n^{-1}S_n= a_n^{-1}\sum_{i=1}^{k_n} (S_{r_n\,i}- S_{r_n\,(i-1)})+o_\P(1)\std \xi_{\alpha}\,.
\eeam
Moreover, since $a_n/a_{r_n}\to\infty$, Theorem~\ref{thm:heyde} yields
\beam\label{eq:ldsum}
\dfrac{\P(\pm a_n^{-1}S_{r_n}>x)}{r_n\,\P(|X|>a_n)}\sim \dfrac{\P(\pm X>x\,a_n)}{\P(|X|>a_n)}\to p_{\pm}\,x^{-\alpha}\,,\qquad x>0\,.  
\eeam
Classical limit theory for triangular arrays of the row-wise 
iid \rv s $(S_{r_ni}-S_{r_n\,(i-1)})_{i=1,\ldots,k_n}$ (e.g. Petrov \cite{petrov:1972}, Theorem 8 
in Chapter IV)
yields that \eqref{eq:sum} holds \fif\  
\beam\label{eq:vague1}
k_n\,\P(a_n^{-1} S_{r_n}\in\cdot )\stv\ \mu_\alpha(\cdot)\,,\\
\lim_{\delta\downarrow 0}\limsup_{\nto}k_n \var\big(a_n^{-1} S_{r_n}\1\big(|S_{r_n}|\le \delta a_n\big)\big)=0\,,  \label{eq:vague2}
\eeam
where $\mu_\alpha$ is defined in \eqref{eq:mualpha}. We notice that \eqref{eq:vague1}
is equivalent to \eqref{eq:ldsum}. 
\par
An alternative way of proving limit theory for the sum process 
$(S_n)$ with an $\alpha$-stable limit $\xi_\alpha$ would be to {\em assume} the 
relations \eqref{eq:vague1} and \eqref{eq:vague2}.
This would be rather indirect and complicated in the case of iid $(X_i)$.
However,  this approach has some merits in the case when $(X_i)$ is a strictly 
stationary \seq\ with a \regvary\ dependence structure, i.e., its \fidi s 
satisfy a multivariate \regvar\ condition (see Davis and Hsing \cite{davis:hsing:1995} or Basrak and Segers \cite{basrak:segers:2009}), and a weak dependence
assumption of the type
\beam\label{eq:mixing}
\E\big[\exp\big(a_n^{-1} it S_n\big)\big]
-\Big(\E\big[\exp\big(a_n^{-1} it S_{r_n}\big)\big]\Big)^{k_n}\to 0\,,\quad 
t\in\bbr\,,\qquad \nto\,,
\eeam
holds. Then $a_n^{-1}S_n\std \xi_\alpha$ \fif\ $a_n^{-1} \sum_{i=1}^{k_n} S_{ni}
\std \xi_\alpha$ where $(S_{ni})_{i=1,\ldots,k_n}$ is an iid \seq\ with the same \ds\ as $S_{r_n}$. Condition \eqref{eq:mixing} is satisfied under 
mild conditions on $(X_i)$, in particular under standard mixing
conditions such as $\alpha$-mixing. Thus one has to
prove the conditions  \eqref{eq:vague1} and \eqref{eq:vague2}. 
In the dependent case the limit \ms\ $\mu_\alpha$ has to be modified: the following analog of \eqref{eq:ldsum} holds: there exists a positive number $\theta_X$ \st
\beao
\dfrac{\P(\pm a_n^{-1}S_{r_n}>x)}{r_n\,\P(|X|>a_n)}\sim \theta_X\,\dfrac{\P(\pm X>x\,a_n)}{\P(|X|>a_n)}\to \theta_X\,p_{\pm}\,x^{-\alpha}\,,\qquad x>0\,.  
\eeao
The quantity $\theta_X$ has an explicit structure in
terms of the so-called tail chain of the \regvary\
\seq\ $(X_i)$. It has interpretation as a {\em cluster index} in the context
of the partial sum operation acting on $(X_i)$.
For details we refer to Mikosch and Wintenberger
\cite{mikosch:wintenberger:2016} and the references therein.
\subsubsection{Extension to the multivariate \regvary\ case.}\label{eq:subsecmultrv}

Consider a \seq\ $(\bfX_i)$  of iid $\bbr^{d}$-valued random vectors with generic element $\bfX$, and define 
\beao
\bfS_0=\boldsymbol{0}\,, \qquad \bfS_n= \bfX_1 + \cdots + \bfX_n\, , \qquad n\ge 1.
\eeao
We say that $\bfX$ is regularly varying with index $\alpha>0$ and a Radon measure $\mu$ on $\bbr^{d}_{\boldsymbol{0}}=\bbr^{d} \backslash \{ {\boldsymbol{0}}\}$, and we write $\mathbf{X} \in \text{RV} (\alpha, \mu )$, if the following vague convergence relation is satisfied on $\bbr^{d}_{\boldsymbol{0}}$:
\beam\label{eq:mrv}
\dfrac{\P(x^{-1} \bfX \in \cdot )}{\P(| \bfX|>x )} \stv \mu(\cdot)\,, \qquad \xto\, ,
\eeam
and $\mu$ has the homogeneity property $\mu(t\, \cdot)=t^{-\alpha}\mu( \cdot)$, $t>0$. We will also use the sequential version of regular variation: for a \seq\ $(a_n)$ \st\ $n\P(| \bfX|>a_n)\to 1$, (\ref{eq:mrv}) is equivalent to
\beao
n\,\P(a_{n}^{-1} \bfX \in \cdot )\stv \mu(\cdot)\,, \qquad \nto\, .
\eeao
For more reading on multivariate \regvar , we refer to Resnick \cite{resnick:1987,resnick:2007}.
\par
Hult et al. \cite{hult:lindskog:mikosch:samorodnitsky:2005}
extended Nagaev's Theorem \ref{thm:nagaev} to the multivariate case: 
\bth[A multivariate Nagaev-type large deviation result]\label{thm:nagaevmv} 
Consider an iid $\bbr^{d}$-valued  \seq\ $(\bfX_i)$ with generic element $\bfX$. Assume the following conditions.
\begin{enumerate}
\item[\rm (1)]
$\mathbf{X} \in \text{\rm RV} (\alpha, \mu )$.
\item[\rm (2)]
The \seq\ of positive numbers $(x_n)$ satisfies
\beam\label{eq:wlln}
	x_{n}^{-1} \bfS_n \stp \boldsymbol{0}\, \mbox{ as } \, \nto\,,
\eeam
and, in addition, 
\beam\label{eq:nagm}
 \left\{\barr{ll}\dfrac{x_{n}^{2}}{n \E[ | \bfX|^{2} \1( | \bfX| \leq x_{n} ) ] \log x_n } \to \infty & \alpha=2 \, \mbox{ and } \, \E[ | \bfX|^{2} ]=\infty\,,\\[4mm]
\dfrac{x_{n}^{2}}{n \log n} \to \infty, & \alpha>2\,\mbox{ or }\, [\alpha=2\,\mbox{ and } \,\E[ | \bfX|^{2} ]<\infty]\,.
\earr\right.
\eeam
\end{enumerate}
 Then 
 \beao
 	\dfrac{\P( x_{n}^{-1} \bfS_n \in \cdot )}{n \P(   |\bfX| >x_{n})} \stv \mu (\cdot )\,,\qquad \nto\,.
\eeao
\ethe

\bre Condition \eqref{eq:wlln} requires that $n\,\E[\bfX]/a_{np}\to\bf0$ for $\alpha>1$. It is always satisfied if $\E[\bfX]=\bf0$. Now assume that the latter condition is satisfied if the expectation of $\bfX$ is finite. If $\alpha\in (0,2)$ 
we can choose any $(p_n)$ \st\ $p_n\to\infty$. 
If $\alpha\ge 2$ and $(np)^{1/\alpha}/n^{0.5+\gamma/\alpha}\to\infty$, equivalently,
$p/n^{\alpha/2-1+\gamma}\to \infty$
 holds for any small $\gamma>0$ then \eqref{eq:nagm} is satisfied.
\ere
The following result extends Theorem \ref{thm:2} to the multivariate case. 
\bth\label{thm:3}
Assume that $\bfX$ satisfies the conditions of Theorem~\ref{thm:nagaevmv}.
Consider an integer \seq\ $p=p_n\to\infty$ and, in addition for $\alpha\ge 2$, 
that $x_n=a_{np}$ satisfies \eqref{eq:nagm}. 
Then the following limit relation holds 
\beam\label{eq:ppconm} 
N _p=\sum_{i=1}^p \vep_{a_{np}^{-1}\bfS_{ni}}\std N\,,
\eeam
where $(\bfS_{ni})$ are iid copies of $\bfS_n$ and $N$ is $\PRM(\mu)$ on $\bbr^d_{\bf0}$.
\ethe
\begin{proof} In view of Proposition~\ref{prop:conv} it suffices to show
that
\beao
p\, \P( a_{np}^{-1} \bfS_n \in \cdot )\stv \mu(\cdot ).
\eeao
Assume $\alpha<2$. Then for 
any \seq\ $p_n\to\infty$, $a_{np}/a_n\to \infty$. Therefore  Theorem
~\ref{thm:nagaevmv} and the definition of $(a_{np})$ imply that for any $\mu$-continuity set $A\subset \bbr^{d}_{\boldsymbol{0}} $,
\beao
p\, \P( a_{np}^{-1} \bfS_n \in A)\sim p\,n\,\P( |\bfX|>a_{np}) \,\mu(A)\to \mu(A)\,,\qquad \nto\,.
\eeao
If $\alpha\ge 2$ the same result holds by virtue of Theorem~\ref{thm:nagaevmv}
and the additional condition \eqref{eq:nagm}. 
\end{proof}
\bexam
Write 
\beao
\bfS_{ni}&=&\big(S_{ni}^{(1)},\ldots,S_{ni}^{(d)}\big)^\top\,,\\
\bfM_n&=& \big(\max_{i=1,\ldots,p}S_{ni}^{(1)},\ldots,\max_{i=1,\ldots,p}S_{ni}^{(d)}\big)^\top= \big(M_n^{(1)},\ldots,M_n^{(d)}\big)^\top\,.
\eeao
For vectors $\bfx,\bfy\in\bbr^d$ with non-negative components, we write
$\bfx\le \bfy$ for the componentwise ordering, $[{\bf0},\bfx]=\{\bfy:{\bf0}\le \bfy\le \bfx\}$ and $[{\bf0},\bfx]^c= \bbr_+^d\backslash [\bf0,\bfx]$. We have
by Theorem~\ref{thm:3},
\beao
\P\big({\bf0}\le a_{np}^{-1}\bfM_n\le\bfx\big)&=& \P\big(N_p([{\bf0},\bfx]^c)=0\big)\\
&\to &\P \big(N([{\bf0},\bfx]^c)=0\big)\\
&=&\exp\big(-\mu([{\bf0},\bfx]^c)\big)=:H(\bfx)\,,\qquad \nto\,,
\eeao
for the continuity points of the \fct\ $-\log H(\bfx)=\mu([{\bf0},\bfx]^c)$.
If $\mu(\bbr_+^d\backslash \{\bf0\})$ is not zero 
$H$ defines a \ds\ on $\bbr_+^d$ with the property $-\log H(t\bfx)=t^{-\alpha}(-\log H(\bfx))$, $t>0$. The non-degenerate components of $H$ are in the type
of the Fr\'echet \ds ; $H$ is referred to as a multivariate Fr\'echet
\ds\ with exponent \ms\ $\mu$.
\eexam

\subsubsection{An extension to iid random sums.}
In this section we consider an alternative random sum process:
\beao
S(t)=\sum_{i=1}^{ \nu(t)}X_{i} \, ,\qquad t\ge 0\,,
\eeao
where $(\nu(t))_{t\ge 0}$ is  a process of integer-valued non-negative 
\rv s independent of the iid \seq\ $(X_i)$ with generic element $X$ and
finite expectation. 
Throughout we assume that $\lambda(t)=\E[\nu(t)]$, $t\ge 0$, is finite 
but $\lim_{t\to\infty}\lambda(t)=\infty$.  We also define
\beao
	m(t)=\E[S(t)]=\E[X]\,\lambda(t) \, .
\eeao
In addition, we assume some technical conditions on the process $\nu$:
\begin{enumerate}
\item[\bf N1]$\nu(t)/\lambda(t) \stp 1$, $t\to\infty$.
\item[\bf N2] There exist $\epsilon, \delta>0$ \st\ 
\beao
\lim_{t\to\infty}\sum_{k>(1+\delta) \lambda(t) } \P (\nu(t)>k)\, (1+\epsilon)^{k} =0 \, .
\eeao
\end{enumerate}
These conditions are satisfied for a wide variety of processes $\nu$, including the  homogeneous Poisson process on $(0,\infty)$.
Kl\"uppelberg and Mikosch \cite{klueppelberg:mikosch:1997} 
proved the following large deviation result for the random sums $S(t)$.
(\cite{klueppelberg:mikosch:1997} allow for the more general condition of extended \regvar .) 
\bth\label{thm:kluppelberg:mikosch}
Assume that $\nu$ satisfies {\rm\bf N1,N2} and is independent of the 
iid non-negative \seq\ $(X_i)$ which is \regvary\ with index $\alpha>1$.
Then for any $ \gamma>0$,
\beao
\sup_{x\ge \gamma\lambda(t)}	\Big|\dfrac{\P(S(t)-m(t) > x)}{\lambda(t) \P(X>x)}-1\Big|\,,\qquad t\to\infty \,.
\eeao
\ethe
The same method of proof as in the previous sections in combination
with the \ld\ result of Theorem~\ref{thm:kluppelberg:mikosch} yields
the following statement. As usual, we assume that $(a(t))$ is a 
\fct\ \st\ $t\,\P(X>a(t))\to 1$ as $t\to\infty$.
\bco\label{thm:6}
Assume the condition of Theorem~\ref{thm:kluppelberg:mikosch}.
Let $(p(t))$ be 
an integer-valued \fct\ \st\ that $p(t)\to\infty$ as $t\to\infty$ 
and a growth condition is satisfied
for every fixed $\gamma>0$ and sufficiently large $t\ge t_0$:
\beam\label{eq:nagr2}
a({\lambda(t) p(t) })\ge \gamma\, \lambda(t)\,.
\eeam 
Then the following limit relation holds for iid copies $S_i$ of the 
random sum process $S$:
\beao
N_{p(t)}=\sum_{i=1}^{p(t)} \vep_{\frac{S_{i}(t)-m(t)}
{a( \lambda(t) p(t))}}\std 
N=\sum_{i=1}^\infty\vep_{\Gamma_i^{-1/\alpha}}\,,\qquad t\to\infty\,,  
\eeao
where $(\Gamma_i)$ is defined in Theorem~\ref{thm:1}.
\eco
\begin{proof} In view of Proposition~\ref{prop:conv}
the result is proved if we can show that as $t\to\infty$,
\beao
p(t) \, \P( (a(\lambda(t) p(t)))^{-1}\big( S(t)- m(t))>x\big)
&\sim & \lambda(t)\,p(t)\,\P( X>a(\lambda(t) p(t) )\,x)
\to  x^{-\alpha}\,,\\
 p(t) \, \P( (a(\lambda(t) p(t)))^{-1}\big( S(t)- m(t))<-x\big)
&\to& 0\,,\qquad x>0\,.
\eeao
But this follows by an application of Theorem~\ref{thm:kluppelberg:mikosch}
in combination with \eqref{eq:nagr2} and the \regvar\ of $X$.
\end{proof}
\bre
Since $a(\lambda(t) p(t))=(\lambda(t)p(t))^{1/\alpha}\ell(\lambda(t)p(t))$ for a 
\slvary\ \fct\ $\ell$ and $\ell(x)\ge x^{-\epsilon/\alpha}$ for any small
$\epsilon>0$ and sufficiently large $x$, \eqref{eq:nagr2} holds
if $p(t)\ge (\lambda(t))^{\alpha-1+\epsilon'}$ for any choice of $\epsilon'>0$ .
\ere

\subsection{An extension: the index of the \pp\ is random}
Let $(P_n)_{n\ge 0}$ be a \seq\ of positive integer-valued \rv s. We assume 
that there exists a \seq\ of positive numbers $(p_n)$ \st\ $p_n\to\infty$ and 
\beam\label{eq:convp}
	\frac{P_n}{p_n} \stp 1 \,,\qquad \nto\,. 
\eeam
This condition is satisfied for wide classes of integer-valued \seq s 
$(P_n)$,
including the renewal counting processes and (inhomogeneous) 
Poisson processes when calculated at the positive integers.
In particular, for renewal processes $p_n\sim c\,n$ provided the
inter-arrival times have finite expectation.
\par
We have the following analog of Proposition~\ref{prop:conv}.
\bpr\label{prop:ranconv}
Let   $(X_{ni})_{n=1,2,\ldots;i=1,2,\ldots}$ be a triangular array of iid \rv s
assuming values in some state space $E\subset\bbr^d$ equipped with the 
Borel $\sigma$-field $\mathcal E$. Let $\mu$ be a Radon \ms\ on ${\mathcal E }$.
If the relation
\beam\label{eq:ranvague}
p_n \,\P\big(X_{n1}\in \cdot\big)\stv \mu(\cdot)\,,\qquad \nto\,,
\eeam
holds on $E$ then
\beao
\wt N_p= \sum_{i=1}^{P_n}\vep_{X_{ni}}\std N\,,\qquad \nto\,,
\eeao
where $N$ is $\PRM(\mu)$ on $E$.  
\epr

\begin{proof} We prove the result by showing \con\ of the Laplace
\fct als. The arguments of a Laplace \fct al are elements of
\beao
C_{K}^{+}(E)= \{ g :E \to \R_{+}: g \mbox{ continuous with compact support} \}\,.
\eeao
For $f\in C_K^+$ we have by independence of the $(X_{ni})$,
\beao
\E\Big[ \exp\Big(-\int_E f\,d \wt N_p\Big)\Big] &=& 
\E\Big[ \exp\Big(-\sum_{j=1}^{P_n}f(X_{nj}) \Big)\Big] = 
\E\Big[ \Big( \E\big[  \exp(-f(X_{n1}) \big] \Big) ^{P_n} \Big] \, .
\eeao
In view of \eqref{eq:convp} there is a real \seq\ $\epsilon_n\downarrow 0$ \st
\beam\label{eq:weak}
\lim_{\nto}\P\big(|P_n/p_n-1|>\epsilon_n\big)=\P(A_n^c)=0\,.
\eeam
Then
\beao\lefteqn{
\E\Big[ \Big( \E\big[  \exp(-f(X_{n1}) \big] \Big) ^{P_n} \Big]}\\
&=& \E\Big[ \Big( \E\big[  \exp(-f(X_{n1}) \big] \Big) ^{P_n} \big(\1(A_n^c)+\1(A_n)\big)
\Big]\\
&=&I_1+I_2\,.
\eeao
By \eqref{eq:weak} we have $I_1\le \P(A_n^c)\to 0$ as $\nto$ while
\beam\label{eq:dominat}
\lefteqn{ \E\Big[ \Big( \E\big[  \exp(-f(X_{n1}) \big] \Big) ^{(1+\epsilon_n)p_n}\1(A_n)\Big]}\nonumber\\&\le&
I_2\le  \E\Big[ \Big( \E\big[  \exp(-f(X_{n1}) \big] \Big) ^{(1-\epsilon_n)p_n}\1(A_n)\Big]\,.
\eeam
In view of Proposition~\ref{prop:conv} and \eqref{eq:ranvague}
\beao
\Big( \E\big[  \exp(-f(X_{n1}) \big] \Big) ^{(1\pm\epsilon_n)p_n}
\to \E\Big(\exp\Big(-\int_E (1-{\rm e}^{-f(\bfx)})\,\mu(d\bfx)\Big)\Big)\,.
\eeao
The \rhs\ is the Laplace \fct al of a $\PRM(\mu)$. 
Now an application of dominated \con\ to $I_2$ in \eqref{eq:dominat} yields
the desired \con\ result. 
\end{proof}
An immediate con\seq\ of this result is that all \pp\ \con s in 
Section~\ref{sec:main} remain valid if the \pp es $N_p$ are replaced by 
their corresponding analogs $\wt N_p$ with a random index 
\seq\ $(P_n)$ independent  of $(S_{ni})$ and satisfying \eqref{eq:convp}. 
Moreover, the growth rates for $p_n\to\infty$ remain the same.
 
\subsection{Extension to the tail empirical process}
We assume that $(S_{ni})$ are iid copies of a real-valued random walk $(S_n)$.
Instead of the \pp es considered in the previous sections one can also
study the tail empirical process
\beao
N_p= \dfrac{1 }{k} \sum_{i=1}^{p} \vep_{c_{[p/k]}^{-1}(S_{ni}/\sqrt{n}-d_{[p/k]})}
\eeao
where $k=k_n\to\infty$, $p=p_n\to\infty$ and $p_n/k_n\to \infty$, and $(c_n)$ and $(d_n)$ are suitable
normalizing and centering constants. To illustrate the theory
we consider two examples.
\bexam 
Assume the conditions and notation of Theorem~\ref{thm:1}. In this case, choose 
$c_n=1/d_n$. Then
\beao
\E[N_p(x,\infty)]&=& \dfrac{p}{k}\,\P\big( S_n/ \sqrt{n}>d_{[p/k]}+x/d_{[p/k]}\big)
\to {\rm e}^{-x}\,,\\
\var\big(N_p(x,\infty)\big)&\le & \dfrac{p}{k^2} \,\P\big( S_n/ \sqrt{n}>d_{[p/k]}+x/d_{[p/k]}\big)\to 0\,,\qquad x\in\bbr, \qquad \nto\,\,,
\eeao
provided $p/k<\exp(\gamma_n^2/2)$.
It is not difficult to see that
\beao
N_p\stp -\log \Lambda\,.
\eeao
Similarly, assume the conditions and the notation of Theorem~\ref{thm:2} and consider 
\beao
N_p= \dfrac{1 }{k} \sum_{i=1}^{p} \vep_{a_{[np/k]}^{-1}(S_{ni}-d_n)}\,.
\eeao
Then for $x>0$ as $\nto$,
\beao
\mathbb{E}[N_p(x,\infty)]&=& \dfrac{p}{k}\,\P\big(a_{[np/k]}^{-1}(S_n-d_n)>x\big)
\sim \dfrac{np}{k}\,\P(X>a_{[np/k]}\, x)\\
&\to& p_+\,x^{-\alpha}=\mu_\alpha(x,\infty)\,,\\
\var(N_p(x,\infty))&\to & 0\,,\\
\mathbb{E}[N_p(-\infty,-x]]&=& \dfrac{p}{k}\,\P\big(a_{[np/k]}^{-1}(S_n-d_n)\le -x\big)
\to p_-\,x^{-\alpha}=\mu_\alpha(-\infty,-x]\,,\\
\var(N_p(-\infty,-x])&\to &0\,,
\eeao
provided the modified \seq\ $p_n/k_n\to\infty$ satisfies the conditions imposed on $(p_n)$
in   Theorem~\ref{thm:2}. We notice that the values of $\mu_\alpha$ on 
$(-\infty,-x]$ and $(x,\infty)$ determine a Radon \ms\ on $\bbr\backslash\{0\}$.
From these relations we conclude that 
 $N_p\stp \mu_\alpha$. Then, following the lines of Resnick and \sta\ \cite{resnick:starica:1995}, Proposition~2.3,
one can for example prove consistency of the Hill estimator based on 
the sample $(S_{ni})_{i=1,\ldots,p}$: assuming for simplicity $d_n=0$, $p_+>0$, we  
write $S_{n,(1)}\ge\cdots\ge  S_{n,(k)}$ for the $k$ largest values. Then
\beao
\dfrac {1}{k} \sum_{i=1}^k \log \dfrac{S_{n,(i)}}{S_{n,(k)}}\stp \dfrac 1 \alpha\,.
\eeao 

\eexam

\subsection{Some related results}\label{subsec:related}
The largest values of sequences of iid  normalized and centered 
partial sum processes  play a role in the context
of random matrix theory which is also the main motivation 
for the present work.  Consider a double array $(X_{it})$ of
iid \regvary\ \rv s with index $\alpha\in (0,4)$ 
(see \eqref{eq:regvar})
and generic element $X$, and also
assume that $\E[X]=0$ if this expectation is finite.
Consider the data matrix
\beao
\bfX:=\bfX_n=(X_{it})_{i=1,\ldots,p;t=1,\ldots,n}
\eeao
and the corresponding sample covariance matrix $\bfX\bfX^\top=(S_{ij})$.
Heiny and Mikosch \cite{heiny:mikosch:2017} proved that
\beam\label{eq:oo}
a_{np}^{-2}\|\bfX\bfX^\top -{\rm diag}(\bfX\bfX^\top)\|_2\stp 0\,,\qquad \nto\,,
\eeam
where $\|A\|_2$ denotes the spectral norm of a $p\times p$ symmetric matrix $A$, 
${\rm diag}(A)$ consists of the diagonal of $A$, $(a_k)$ is any \seq\
satisfying $k\,\P(|X|>a_k)\to 1$ as $k\to\infty$, and $p_n=n^\beta\ell(n)$ for
some $\beta\in (0,1]$ and a \slvary\ \fct\ $\ell$.
Write $\la_{(1)}(A)\ge \cdots\ge \la_{(p)}(A)$ for the ordered eigenvalues of $A$.
According to Weyl's inequality (see Bhatia \cite{bhatia:1997}), the eigenvalues of $\bfX\bfX^\top$
satisfy the relation
\beam\label{eq:weyl}
a_{np}^{-2}\sup_{i=1,\ldots,p}\big|\la_{(i)}(\bfX\bfX^\top)- \la_{(i)}(\diag(\bfX\bfX^\top)) \big|\le a_{np}^{-2}\|\bfX\bfX^\top -{\rm diag}(\bfX\bfX^\top)\|_2\stp 0\,.
\eeam
But of course, $\la_{(i)}(\diag(\bfX\bfX^\top))$ are the ordered values of
the iid partial sums $S_{ii}=\sum_{t=1}^n X_{it}^2$, $i=1,\ldots,p$. In view
of \eqref{eq:weyl} the \asy\ theory for the largest eigenvalues of the 
normalized sample covariance matrix $a_{np}^{-2} \bfX\bfX^\top$
(which also needs centering for $\alpha\in (2,4)$)  are determined
through the Fr\'echet \con\ of the processes with points 
$(a_{np}^{-2}S_{ii})_{i=1,\ldots,p}$.
Moreover, \eqref{eq:weyl} implies the Fr\'echet \con\ of 
the \pp es of the normalized and centered eigenvalues
of the sample covariance matrix.
\par
\par
The \ld\ approach also works for proving limit theory for the \pp\ of the 
off-diagonal elements of $\bfX\bfX^\top$ provided $X$ has sufficiently high 
moments.
 Heiny et al. \cite{heiny:mikosch:yslas:2019} prove Gumbel \con\ 
for the \pp\ of the off-diagonal elements $(S_{ij})_{1\le i<j\le p}$. 
The situation is more complicated because the points $S_{ij}$
are typically dependent. Multivariate extensions of the  normal \ld\ approximation $0.5p^2 \P(d_{p^2/2}(S_{12}-d_{p^2/2})>x)\to \exp(-x)$ show that
the \pp\ of the standardized $(S_{ij})$ has the same limit Poisson process 
as if the $S_{ij}$ were independent. Moreover, \cite{heiny:mikosch:yslas:2019}
show that the \pp\ of the diagonal elements $(S_{ii})$ (under suitable 
conditions on the rate of $p_n\to\infty$ and under $\E[|X|^s]<\infty$ for $s>4$) 
converges to $\PRM(-\log \Lambda)$.
This result indicates that the off-diagonal and diagonal entries
of $\bfX\bfX^\top$ exhibit very similar extremal behavior.
This is in stark contrast to the aforementioned results in 
 \cite{heiny:mikosch:2017} where the diagonal entries have Fr\'echet
extremal behavior.
\par
Related results can also be found 
in Gantert and H\"ofelsauer \cite{gantert:hoefelsauer:2018} who
consider real-valued branching random walks and prove a large deviation principle for the position of the right-most particle;  see Theorem 3.2 in \cite{gantert:hoefelsauer:2018}.
The position of the right-most
particle is the maximum of a collection of a random number of dependent random
walks. In this context, the authors also prove a related \ld\ result under the assumption that the considered random walks are iid. They show that 
the maximum of these iid random walks stochastically dominates the maximum of 
the branching random walks; see Theorem 3.1 and Lemma 5.2 in \cite{gantert:hoefelsauer:2018}.  An early comparison between maxima of branching and iid 
random walks was provided by Durrett \cite{durrett:1979}.

\end{document}